\documentclass[a4paper]{article} 

\usepackage{latexsym}
\usepackage{amsfonts}
\usepackage{amssymb}
\usepackage{amsmath}
\usepackage{amscd}
\usepackage{eucal}
\usepackage[all, knot]{xy}
\xyoption{arc}
\xyoption{rotate}
\xyoption{dvips}
\usepackage{amsthm}
\usepackage[dvips]{pict2e}

\usepackage{verbatim} 
\usepackage[dvips]{graphicx}

\newtheorem{Def}{Definition}[section]

\newtheorem{prop}[Def]{Proposition}

\newtheoremstyle{example}{\topsep}{\topsep}%
     {}
     {}
     {\bfseries}
     {.}
     {2pt}
     {\thmname{#1}\thmnumber{ #2}\thmnote{ #3}}

   \theoremstyle{example}

\newcommand{\lra}{\ensuremath{\longrightarrow}}
\newcommand{\cat}[1]{\ensuremath{\mbox{\bfseries {\upshape {#1}}}}}
\newcommand{\cl}[1]{\ensuremath{\mathcal {#1}}}
\newcommand{\bb}[1]{\ensuremath{\mathbb {#1}}}
\newcommand{\ed}{\end{document}}
\newcommand{\map}[1]{\ensuremath{\stackrel{{#1}}{\lra}}}

\newcommand{\Set}{{\cat{Set}}}

\newcommand{\id}{\ensuremath{\mbox{\em id}}}

\newcommand{\numarabic}{\renewcommand{\labelenumi}{\arabic{enumi}.}}

\newcommand{\sunit}{\setlength{\unitlength}{1mm}}
\newcommand{\setunit}[1]{\setlength{\unitlength}{#1}}

\newenvironment{prf}{\vspace{2ex}\begin{sloppypar}{\noindent\upshape
{\bfseries Proof. }}} {{\hspace*{\fill}
$\Box$}\end{sloppypar}\vspace{2ex}}

\newcommand{\adj}{\dashv}

\newcommand{\rrr}{
\begin{picture}(60,70)(0,15)

\qbezier(2,9)(26,32)(35,60)
\qbezier(35,60)(36,62)(35.6,64)
\qbezier(35.6,64)(35.2,68.2)(32.2,65)
\qbezier(32.2,65)(29,61.5)(28.5,55)
\qbezier(28.5,55)(27.4,42.4)(41,46)
\qbezier(41,46)(50.2,48.8)(58.4,57)
\qbezier(58.4,57)(51.6,40)(52,28)
\qbezier(52,28)(52.6,12.2)(64.8,9.5)
\qbezier(64.8,9.5)(71.5,7.8)(77,14)
\end{picture}
}

\newcommand{\curlyr}{\setunit{0.03mm}
\makebox[2.8mm][l]{\begin{picture}(110,70)(38,-10)
\put(0,0){\rrr}
\put(1,0){\rrr}
\put(2,0){\rrr}
\end{picture}}}

\newcommand{\llll}{
\begin{picture}(65,80)
\qbezier(5,7)(9,4)(15,6)
\qbezier(15,6)(23.5,8.5)(29,15)
\qbezier(29,15)(38,26)(42.8,40)
\qbezier(42.8,40)(46,49)(47,60)
\qbezier(47,60)(48.5,65.6)(47.8,71)
\qbezier(47.8,71)(47.5,77.5)(44,72)
\qbezier(44,72)(41.5,68)(40.5,60)
\qbezier(40.5,60)(37.5,41)(38.5,30)
\qbezier(38.5,30)(38.5,22)(40,18)
\qbezier(40,18)(42,12)(46,8)
\qbezier(46,8)(52,1.5)(58.5,8)
\end{picture}}

\newcommand{\curlyl}{\setunit{0.04mm}
\makebox[2.9mm][l]{\begin{picture}(80,70)(28,5)
\put(0,0){\llll}
\put(1,0){\llll}
\put(2,0){\llll}
\end{picture}}}

\newcommand{\bdot}[1]
{\sunit
\begin{picture}(0,0)
\put(0,0){\circle*{#1}}
\end{picture}}

\newcommand{\bigleft}{\left(\begin{array}{l}\vspace{3mm}\end{array}\right.
\hspace{-5mm}}

\newcommand{\bigright}{\hspace{-5mm}\left.\begin{array}{l}\vspace{3mm}
\end{array}\right)}

\newcommand{\op}{\ensuremath{\mbox{\hspace{1pt}{\scriptsize\upshape op}}}}

\newcommand{\tcurve}{\qbezier(0,0)(20,35)(40,0)}
\newcommand{\bcurve}{\qbezier(0,0)(20,-35)(40,0)}
\newcommand{\flatbit}{\line(1,0){40}}
\newcommand{\tlabel}[1]{\put(20,6){\makebox(0,0){#1}}}
\newcommand{\blabel}[1]{\put(20,-4){\makebox(0,0)[t]{#1}}}
\newcommand{\mlabel}[1]{\put(20,0){\makebox(0,0){#1}}}

\newcommand{\twelveoc}[2]{
\put(0,0){\tcurve}
\put(0,0){\bcurve}
\put(0,0){\mlabel{#1}}   
\put(40,0){\tcurve}
\put(40,0){\bcurve}
\put(40,0){\mlabel{#2}}}  

\newcommand{\tcell}[1]{
\put(0,0){\tcurve}
\put(0,0){\flatbit}
\put(0,0){\tlabel{#1}}}

\newcommand{\bcell}[1]{
\put(0,0){\bcurve}
\put(0,0){\flatbit}
\put(0,0){\blabel{#1}}}

\newcommand{\threeoc}[2]{
\put(0,0){\tcell{#1}}
\put(0,0){\bcell{#2}}}

\newcommand{\oneoc}[4]{
\put(-5,0){\threeoc{#1}{#3}}
\put(45,0){\threeoc{#2}{#4}}}

\newcommand{\ttwo}[2]{
\put(0,0){\tcell{#1}}
\put(40,0){\tcell{#2}}}

\newcommand{\btwo}[2]{
\put(0,0){\bcell{#1}}
\put(40,0){\bcell{#2}}}

\newcommand{\twooc}[4]{
\put(0,5){\ttwo{#1}{#2}}
\put(0,-5){\btwo{#3}{#4}}}

\newcommand{\ttl}{\ensuremath{\triangleleft}}
\newcommand{\ttr}{\ensuremath{\triangleright}}
\newcommand{\ttu}{\ensuremath{\scriptstyle\bigtriangleup}}
\newcommand{\ttd}{\ensuremath{\scriptstyle\bigtriangledown}}

\newcommand{\nowdot}{\bdot{1}}

\begin{document}

\title{A note on the Penon definition of $n$-category}
\author{Eugenia Cheng \\Universit\'e de Nice-Sophia Antipolis and University of Sheffield\\and\\Michael Makkai\\ McGill University}

\date{June 2008}
\maketitle


\begin{abstract}

Nous \'etudions les tricat\'egories de Penon, et nous d\'emontrons que, dans le cas des tricat\'egories deux fois d\'eg\'en\'er\'ees, on obtient les cat\'egories mono\"idales sym\'etriques et non les cat\'egories mono\"idales tress\'ees.  Nous prouvons que les tricat\'egories de Penon ne peuvent pas donner toutes les tricat\'egories.  Pour
corriger cette situation, nous proposons une petite modification de la d'efinition, utilisant les ensembles globulaires non-r\'eflexifs \`a la place des ensembles globulaires r\'eflexifs, et nous d\'emontrons qu'ainsi le probl\`eme pr\'ec\'edent relatif aux tricat\'egories deux fois d\'eg\'en\'er\'ees n'appara\^it plus.\\

We show that doubly degenerate Penon tricategories give symmetric rather than braided monoidal categories.  We prove that Penon tricategories cannot give all tricategories, but we show that a slightly modified version of the definition rectifies the situation.  We give the modified definition, using non-reflexive rather than reflexive globular sets, and show that the problem with doubly degenerate tricategories does not arise.\\

\noindent \sloppy {\bfseries Keywords:} tricategory, degenerate tricategory, braided monoidal category, symmetric monoidal category, globular set, reflexive, non-reflexive.\\


\end{abstract}


\section*{Introduction}

Many different definitions of weak $n$-category have been proposed but as yet the relationship between them and their validity have not been well understood.  One preliminary check that can be applied to any proposed definition is that it ``agrees", in some suitable sense, with well-established low-dimensional examples.  Thus, one might begin by checking the definition of 1-category against the usual definition of category (and this is not always trivial) and then the definition of 2-category against the classical definition of bicategory \cite{ben1}.

 After this point things become more difficult; the definition of tricategory \cite{gps1} has been accepted but a completely algebraic version of the definition has since been proposed \cite{gur2}, and questions remain about what should be ``the'' definition of tricategory, if indeed a unique definition should be sought at all.  However, there is a degenerate form of tricategory which is much better understood -- that is, braided monoidal categories.

Corollary 8.7 of \cite{gps1} states that ``one-object, one-arrow tricategories are precisely braided monoidal categories".  However, the results of \cite{cg3} show that the correspondence is not straightforward, and one-object one-arrow tricategories (``doubly degenerate tricategories'') in fact give rise to braided monoidal categories with various extra pieces of structure.  However, it is shown that braided monoidal categories should at least \emph{arise among} the totality of tricategories; by focussing on this in the present work we avoid the intricate questions involved in the above Corollary.

The main aim of this paper is to show that Penon's definition of $n$-category in its original form is not as general as it might be, as it gives symmetric rather than braided monoidal categories in the ``doubly degenerate" case.  This should not be seen as a serious problem -- the situation is quite easily rectified by starting with globular sets instead of reflexive globular sets in Penon's definition.  A reflexive globular set is one in which putative identities are already picked out. This can be thought of as being analogous to degeneracies being part of the structure of a simplicial set, but then the analogy breaks down.  For simplicial sets the presence of degeneracies is a rich and crucial part of the structure, but in Penon's definition of $n$-category, it is precisely these degeneracies that cause the resulting definition to be slightly too strict, yielding a symmetry instead of a braiding in the one-object, one-arrow situation.

The problem in the reflexive case is that, since the identities are picked out in the underlying globular set, when forming the free 3-cateogry on such a structure it is possible to have non-identity cells whose source and target are the identity.  In the non-reflexive case, the identities are added in \emph{freely} in the free 3-category, so the only cells whose source or target are identities are themselves identities.  It is a general principle that having non-trivial cells with identity source and target causes problems, as in the following situations:

\begin{itemize}

\item Strict computads do not form a presheaf category but strict many-to-one computads do; the problem is caused by the possibility of cells with source and target the identity \cite{mz1}, so this is avoided by insisting that the target is 1-ary (thus disallowing identities since they are ``nullary'').

\item Coherence for tricategories \cite{gur1} says that all diagrams of constraints in a \emph{free} tricategory commute, but not all diagrams of constraints in a general tricategory commute; in a general tricategory a diagram of constraints commutes if it involves no non-identity cells with the identity in the source or target.  For example, the diagram asserting that a braiding is in fact a symmetry does not necessarily commute; it involves cells whose source and target are the identity.

\end{itemize}

We present the result in two different ways.  The first, more intuitively clear but less precise, says that ``A degenerate Penon tricategory is a symmetric monoidal category".  This is the subject of Section~\ref{two}.  We simply examine a degenerate Penon tricategory and express it as a braided monoidal category in the expected way; we then see that in fact the braiding is forced to be a symmetry.  However, this is not a precise mathematical statement -- all it says is that the generally expected method of producing a braided monoidal category from a doubly degenerate tricategory does not gives us \emph{all} braided monoidal categories, only the symmetric ones.  However this may be considered to be the heart of the problem, and was pointed out by the second author during the Workshop on $n$-categories at the IMA in June 2004.

In Section~\ref{three} we ``go backwards" in order to make a precise statement.  First we exhibit a monoidal category which can be equipped with a braiding but \emph{cannot} be equipped with a symmetry.  We then express this as a tricategory and show that it does not satisfy the axioms for a Penon tricategory.  Thus we conclude that Penon's original notion of tricategory does not include all the examples we would like.


In Section~\ref{four} we give the non-reflexive version of Penon's definition, and in Section~\ref{five} we show that these problems do not arise in this case.  This non-reflexive version is essentially that given in \cite{bat3}, although in that work it is conjectured that the reflexive and non-reflexive versions are equivalent.

We begin in Section~\ref{one} by reviewing the basic definitions.  Note that we will often use the term ``$n$-category'' even when $n$ might be $\omega$.

\section{Basic definitions}\label{one}

In this section we recall the definition of $n$-category proposed by Penon \cite{pen1}.  According to this definition, an $\omega$-category is an algebra for a certain monad $P$ on the category of reflexive globular sets.
Our statement of the definition is more similar to that of Leinster \cite{lei7}; for more explanation we also refer the reader to \cite{cl1}.  The definition starts with the underlying data given by a \emph{reflexive globular set}, then imposing the structure of a \emph{magma} (for composition) and \emph{contraction} (for coherence).  For finite $n$ a simple truncation is applied to the underlying data, while some care must be taken over the $n$-cells when defining contractions in this case.

\subsection{Reflexive globular sets}

We write \cat{RGSet} for the category of reflexive globular sets.  \cat{RGSet} is the category of presheaves $[\bb{R}^{\op}, \Set]$ where \bb{R} is the category whose objects are the natural numbers and whose morphisms are as depicted below:

\[\xy
(-20,10)*+{\cdots}="b4";
(0,10)*++{3}="b3";
(20,10)*++{2}="b2";
(40,10)*++{1}="b1";
(60,10)*++{0}="b0";
{\ar@<.5ex>^{t} "b3"; "b4"};
{\ar@<-.5ex>_{s} "b3"; "b4"};
{\ar@<.5ex>^{t} "b2"; "b3"};
{\ar@<-.5ex>_{s} "b2"; "b3"};
{\ar@<.5ex>^{t} "b1"; "b2"};
{\ar@<-.5ex>_{s} "b1"; "b2"};
{\ar@<.5ex>^{t} "b0"; "b1"};
{\ar@<-.5ex>_{s} "b0"; "b1"};
{\ar@/_1.7pc/_{i} "b1";"b0"};
{\ar@/_1.7pc/_{i} "b2";"b1"};
{\ar@/_1.7pc/_{i} "b3";"b2"};

\endxy\]
satisfying globularity and reflexivity conditions.  However we will write a reflexive globular set explicitly as a diagram of sets

\[\xy
(-20,10)*+{\cdots}="b4";
(0,10)*+{A(3)}="b3";
(20,10)*+{A(2)}="b2";
(40,10)*+{A(1)}="b1";
(60,10)*+{A(0)}="b0";
{\ar@<.5ex>^{s} "b4"; "b3"};
{\ar@<-.5ex>_{t} "b4"; "b3"};
{\ar@<.5ex>^{s} "b3"; "b2"};
{\ar@<-.5ex>_{t} "b3"; "b2"};
{\ar@<.5ex>^{s} "b2"; "b1"};
{\ar@<-.5ex>_{t} "b2"; "b1"};
{\ar@<.5ex>^{s} "b1"; "b0"};
{\ar@<-.5ex>_{t} "b1"; "b0"};
{\ar@/^1.7pc/^{i} "b0";"b1"};
{\ar@/^1.7pc/^{i} "b1";"b2"};
{\ar@/^1.7pc/^{i} "b2";"b3"};

\endxy\]
satisfying the globularity conditions $ss = st$ and  $ts = tt$, and the reflexivity condition $si = ti = 1$.  For the finite $n$-dimensional case we truncate the diagram to get the category $\cat{RGSet}_n$ of $n$-dimensional reflexive globular sets as below

\[\xy
(-33,10)*+{A(n)}="b5";
(-10,10)*+{A(n-1)}="b4";
(12,10)*+{\ }="b3";
(16,10)*{\cdots};
(20,10)*+{\ }="b2";
(40,10)*+{A(1)}="b1";
(60,10)*+{A(0)}="b0";
{\ar@<.5ex>^<<<<<<{s} "b5"; "b4"};
{\ar@<-.5ex>_<<<<<<{t} "b5"; "b4"};
{\ar@<.5ex>^>>>>>>>>{s} "b4"; "b3"};
{\ar@<-.5ex>_>>>>>>>>{t} "b4"; "b3"};
{\ar@<.5ex>^<<<<<<<{s} "b2"; "b1"};
{\ar@<-.5ex>_<<<<<<<{t} "b2"; "b1"};
{\ar@<.5ex>^{s} "b1"; "b0"};
{\ar@<-.5ex>_{t} "b1"; "b0"};
{\ar@/^1.7pc/^{i} "b0";"b1"};
{\ar@/^1.7pc/^{i} "b1";"b2"};
{\ar@/^1.7pc/^{i} "b3";"b4"};
{\ar@/^1.7pc/^{i} "b4";"b5"};

\endxy\]

We call the elements of $A(k)$ the $k$-cells of $A$.  The maps $s$ and $t$ give the source and target of each $k$-cell and the map $i$ picks out the putative identity for each $k$-cell.  Part of the structure of the monad $P$ will be to ensure that these really do act as (weak) identities in the $n$-category structure.

A map of reflexive globular sets is a map of these diagrams making all the obvious squares commute.

\paragraph{Note} Every strict $n$-category has an underlying $n$-dimensional reflexive globular set.

\subsection{Magmas}

A magma is a reflexive globular set equipped with binary composition at all levels.   That is, for all $m \geq 1$ we can compose along bounding $k$-cells for any $0 \leq k \leq m-1$.  So given $\alpha, \beta \in A(m)$ with
	\[t^{m-k}\alpha = s^{m-k}\beta\]
we have composite $\beta \circ_k \alpha \in A(m)$ with source and target given by
	\[ s(\beta \circ_k \alpha) = \left\{\begin{array}{cl}
	s(\beta) \circ_k s(\alpha) & \mbox{\hspace{2em}if\ \ } k < m-1 \\
	s(\alpha) & \mbox{\hspace{2em}if\ \ } k = m-1
	\end{array}\right.\]

\[ t(\beta \circ_k \alpha) = \left\{\begin{array}{cl}
	t(\beta) \circ_k t(\alpha) & \mbox{\hspace{2em}if\ \ } k < m-1 \\
	t(\beta) & \mbox{\hspace{2em}if\ \ } k = m-1
	\end{array}\right.\]
Note that the composites on the right hand side make sense because of the globularity conditions. For examples and diagrams illustrating these composites see \cite{cl1}.

A map of magmas is a map of the underlying reflexive globular sets preserving composition.

An $n$-dimensional magma is one whose underlying reflexive globular set is $n$-dimensional.

\paragraph{Note} In a magma only binary composites are given (i.e. in the language of Leinster it is ``biased") and no axioms are required to be satisfied.  In particular, the putative identities are still not required to act as identities with respect to the composition.  Further, note that any strict $n$-category has an underlying $n$-dimensional magma.

\subsection{Contractions}

Composition in a magma is not required to be in any way coherent; we achieve coherence for $n$-categories by way of a ``contraction''.  A contraction is a piece of structure that can be defined on any map $A \map{f} B$ of reflexive globular sets.  The idea of a contraction is similar to contractibility of topological spaces, in that it measures holes, or rather lack thereof.  A contraction on a map $A \map{f} B$ essentially ensures that $A$ has ``no more holes up to homotopy'' than $B$; the contraction cells witness the contraction of $A$ onto $B$.

First we need a notion of parallel $k$-cells.

\begin{Def} A pair of $k$-cells $\alpha, \beta$ are called \emph{parallel} if
	\begin{itemize}
	\item $k = 0$, or
	\item $k \geq 1$ and $s\alpha = s\beta, t\alpha = t\beta$.
	\end{itemize}
\end{Def}

\begin{Def} A \emph{contraction} [\ ,\ ] on a map $A \map{f} B$ of reflexive globular sets gives, for any pair of parallel $k$-cells $\alpha$ and $\beta$ such that $f\alpha = f\beta$, a $(k+1)$-cell
	\[ [\alpha, \beta] : \alpha \lra \beta\]
such that
	\begin{enumerate}
	\item $f[\alpha, \beta] = i(f(\alpha))$
	\item $[\alpha,\alpha] = i\alpha$.
	\end{enumerate}
For the $n$-dimensional case a contraction gives the above for $k < n$; given $\alpha$ and $\beta$ as above for $k=n$ we must have $\alpha = \beta$.
\end{Def}

The cells $[\alpha,\beta]$ are referred to generally as ``contraction cells".  This definition can also be thought of informally as saying that ``any disc in $B$ with a lift of its boundary to $A$ gives a lift of the disc as well''.

\subsection{The crucial category \cl{Q}}

We construct the monad $P$ from an adjunction

\[\xy
(0,20)*++{\cl{Q}}="t0";
(0,0)*++{\cat{RGSet}}="b0";
{\ar@<1.2ex>^{U} "t0"; "b0"};
{\ar@<1.2ex>_{\adj}^{F} "b0"; "t0"};

\endxy\]
The category \cl{Q} has objects of the form

\sunit
\begin{center}
\begin{picture}(30,30)
\put(5,25){$A$}
\put(5,5){$B$}
\put(6,23){\vector(0,-1){14}}
\put(10,15){$f,\ [\ ,\ ]$}
\end{picture}
\end{center}
where

\begin{itemize}
\item $A$ is a magma
\item $B$ is a strict $n$-category
\item $f$ is a map of magmas
\item $[\ ,\ ]$ is a contraction on $f$
\end{itemize}
The idea is that the contraction ensures that $A$ has ``no more holes up to homotopy'' than the strict $n$-category $B$, so although it need not be a strict $n$-category itself, it cannot be too incoherent.

A morphism of such objects is a square

\sunit
\begin{center}
\begin{picture}(30,40)
\put(5,25){$A$}
\put(5,5){$B$}
\put(6,23){\vector(0,-1){14}}
\put(3,15){$f$}

\put(25,25){$A'$}
\put(25,5){$B'$}
\put(26,23){\vector(0,-1){14}}
\put(27,15){$f'$}

\put(9,26){\vector(1,0){15}}
\put(9,6){\vector(1,0){15}}

\put(15,28){$g_1$}
\put(15,3){$g_2$}

\end{picture}
\end{center}

\noindent where $g_1$ is a map of magmas, $g_2$ is a map of strict $n$-categories, and the square of underlying magma maps commutes; furthermore the maps must preserve contraction cells, that is, for every contraction cell $[\alpha,\beta]$ in $A$, we must have
\[ g_1[\alpha,\beta]=[g_1 \alpha, g_1 \beta] .\]

\subsection{The adjunction}

There is a forgetful functor
	\[\cl{Q} \map{U} \cat{RGSet}\]
sending an object as above to the underlying reflexive globular set of $A$.  This functor has a left adjoint $F$; we define $P$ to be the monad induced by this adjunction.  For the $n$-dimensional case we write $P_n$ for the induced monad on $\cat{RGSet}_n$

The existence of this adjoint can be deduced from standard Adjoint Functor Theorems; it is proved by Penon in \cite[Section 4]{pen1}, and it is also quite straightforward to construct using results of \cite{che14}.  Given a reflexive globular set $A$ we can express $FA$ as

\sunit
\begin{center}
\begin{picture}(30,30)
\put(5,25){$PA$}
\put(5,5){$T_RA$}
\put(8,23){\vector(0,-1){14}}
\put(10,15){$\phi,\ [\ ,\ ]$}
\end{picture}
\end{center}
where $T_R$ is the free strict $\omega$-category monad on reflexive globular sets.

\paragraph{Note} This is quite different from a free strict $\omega$-category on a \emph{non-reflexive} globular set; for example $T_R1=1$ as reflexive globular sets (as the unique cell at each dimension must be the identity) which is certainly not true of non-reflexive globular sets.  This difference may be thought of as being the heart of the problem considered in this paper.

\bigskip

The idea is to combine two types of structure: contraction and magma.  We proceed dimension by dimension -- at each level we first add in the required contraction cells freely, and then binary composites freely.  $\phi$ then acts by sending all contraction cells to the identity in $T_RA$, and forgetting the parentheses in all composites.

For the finite $n$-dimensional case the final stage of the construction consists of \emph{identifying} any $n$-dimensional composites that lie over the same cell in $T_RA$.

\begin{Def} An \emph{$\omega$-category} is defined to be a $P$-algebra.  An \emph{$n$-category} is defined to be a $P_n$-algebra.
\end{Def}

\section{Doubly degenerate 3-categories as symmetric monoidal categories}\label{two}

In this section we show how a doubly degenerate Penon 3-category gives rise to a braided monoidal category, and that the braiding given in this way is in fact necessarily a symmetry.  Since the main aim of this section is to show why the braiding must be a symmetry, we do not go through the details of checking all the axioms for a braided monoidal category.

A doubly degenerate tricategory is one that has only one 0-cell and one 1-cell.  The general idea is as follows.  We obtain a category from it by regarding the old 2-cells as objects and the old 3-cells as morphisms, that is, we take the unique hom-category on the unique 1-cell.  We obtain a monoidal structure by taking the tensor product to be given by vertical composition of 2-cells.  Finally we use an ``up to isomorphism" Eckmann-Hilton argument to show that this tensor is ``commutative up to isomorphism" -- that is, it is a braiding.

We now state this in the framework of Penon's definition.  For convenience we now write $P=P_3$ for the ``free Penon 3-category" monad on reflexive 3-globular sets.  Let
$\left(
\def\objectstyle{\scriptstyle}
\def\labelstyle{\scriptstyle}
\vcenter{\xymatrix @-1.0pc  { PA \ar[d]^{\theta} \\ A  }} \right)$
be a $P$-algebra where $A$ is a doubly degenerate reflexive 3-globular set i.e. it has only one 0-cell and only one 1-cell.  We construct a braided monoidal category from it as follows:
	\begin{itemize}
	\item the objects are given by $A(2)$
	\item the morphisms are given by $A(3)$
	\item the tensor product is given by $\alpha \otimes \beta = \alpha \circ \beta$ as 2-cells of $A$
	\item the braiding $\gamma_{\alpha, \beta} : \alpha \otimes \beta \lra \beta \otimes \alpha$ is given by the contraction cell $[\alpha \circ \beta, \beta \circ \alpha]$.
	\end{itemize}
To see that this contraction cell exists we need to show that
	\[\phi(\alpha \circ \beta) = \phi(\beta \circ \alpha) \in T_RA\]
where $\phi$ is the map $PA \lra T_RA$.  So we need to show that $\alpha \circ \beta = \beta \circ \alpha \in T_RA$.  This is proved by an Eckmann-Hilton type argument using the fact that the source and target 1-cells of $\alpha$ and $\beta$ are the identity in $T_RA$.  We find it helpful to place the various stages of the argument on the following ``clock face":

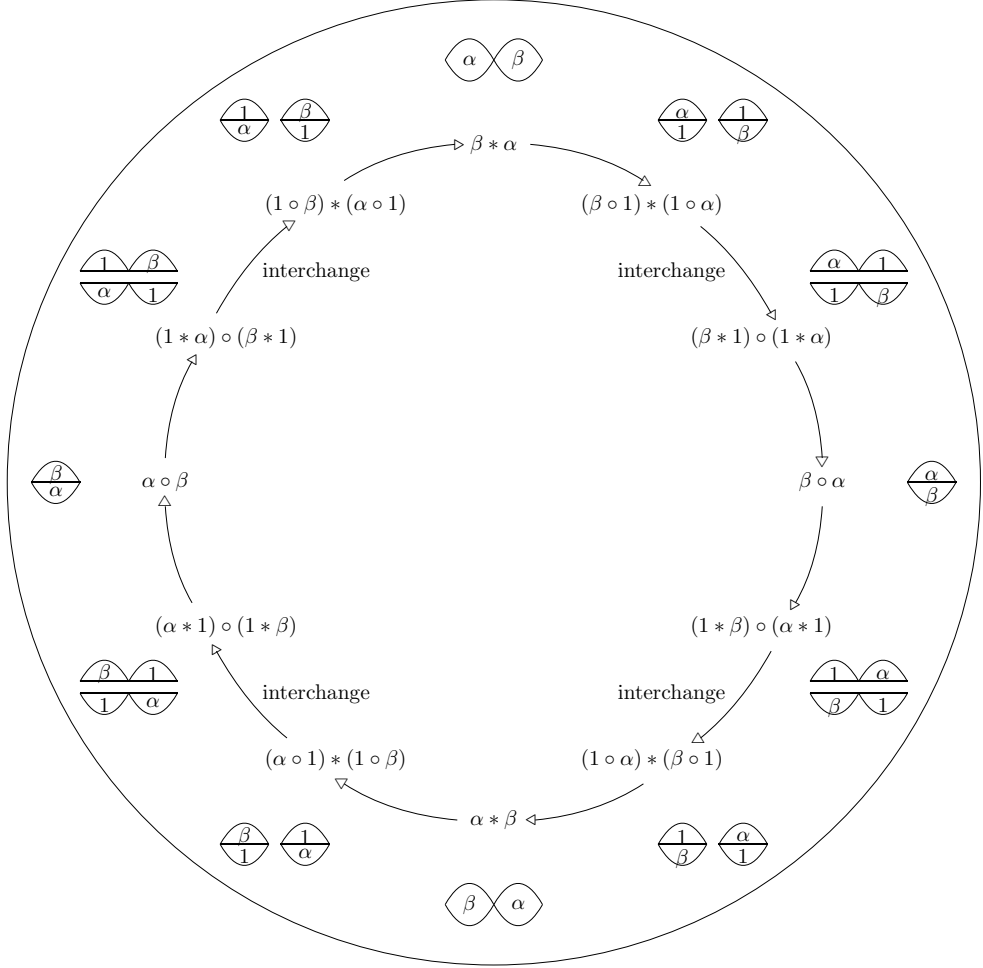
\scalebox{0.8}{\setunit{0.2mm}

\hspace*{10mm}\begin{picture}(800,800)(100,0)

\put(400,400){\circle{800}}

\put(360,750){\twelveoc{$\alpha$}{$\beta$}}
\put(540,700){\oneoc{$\alpha$}{1}{1}{$\beta$}}
\put(660,570){\twooc{$\alpha$}{1}{1}{$\beta$}}
\put(740,400){\threeoc{$\alpha$}{$\beta$}}

\put(660,230){\twooc{1}{$\alpha$}{$\beta$}{1}}
\put(540,100){\oneoc{1}{$\alpha$}{$\beta$}{1}}
\put(360,50){\twelveoc{$\beta$}{$\alpha$}}

\put(180,100){\oneoc{$\beta$}{1}{1}{$\alpha$}}
\put(60,230){\twooc{$\beta$}{1}{1}{$\alpha$}}
\put(20,400){\threeoc{$\beta$}{$\alpha$}}

\put(60,570){\twooc{1}{$\beta$}{$\alpha$}{1}}
\put(180,700){\oneoc{1}{$\beta$}{$\alpha$}{1}}

\put(400,680){\makebox(0,0){$\beta \ast \alpha $}}
\qbezier(430,680)(485,675)(523,650)
\put(518,645.7){\ttu}
\put(530,630){\makebox(0,0){$(\beta \circ 1) \ast (1 \circ \alpha) $}}
\qbezier(570,612)(603,585)(628,540)
\put(623.5,534.4){\ttl}
\put(620,520){\makebox(0,0){$(\beta \ast 1) \circ (1 \ast \alpha) $}}
\qbezier(648,500)(667,467)(670,420)
\put(663.5,415){\ttd}
\put(670,400){\makebox(0,0){$\beta \circ \alpha $}}
\qbezier(648,300)(667,333)(670,380)
\put(643.5,294.5){\ttr}
\put(620,280){\makebox(0,0){$(1 \ast \beta) \circ (\alpha \ast 1)$}}
\qbezier(570,188)(603,215)(628,260)
\put(561.5,183.5){\ttu}
\put(530,170){\makebox(0,0){$(1 \circ \alpha) \ast (\beta \circ 1)$}}
\qbezier(434,120)(485,125)(523,150)
\put(426,116){\ttl}
\put(400,120){\makebox(0,0){$\alpha \ast \beta$}}
\qbezier(370,120)(315,125)(277,150)
\put(268.5,148.5){\ttd}
\put(270,170){\makebox(0,0){$(\alpha \circ 1) \ast (1 \circ \beta)$}}
\qbezier(230,188)(197,215)(172,260)
\put(167.8,257.3){\ttr}
\put(180,280){\makebox(0,0){$(\alpha \ast 1) \circ (1 \ast \beta)$}}
\qbezier(152,300)(133,333)(130,380)
\put(123.5,380){\ttu}
\put(130,400){\makebox(0,0){$\alpha \circ \beta $}}
\qbezier(152,500)(133,467)(130,420)
\put(147.3,497.1){\ttl}
\put(180,520){\makebox(0,0){$(1\ast \alpha) \circ (\beta \ast 1)$}}
\qbezier(230,612)(197,585)(172,540)
\put(225.5,610.5){\ttd}
\put(270,630){\makebox(0,0){$(1\circ \beta) \ast (\alpha \circ 1)$}}
\qbezier(368,680)(315,675)(277,650)
\put(367,675){\ttr}

\put(590,575){\makebox(0,0)[r]{interchange}}
\put(590,225){\makebox(0,0)[r]{interchange}}
\put(210,225){\makebox(0,0)[l]{interchange}}
\put(210,575){\makebox(0,0)[l]{interchange}}


\end{picture}}

\bigskip

Since all these composites are equal in the strict 3-category $T_RA$, we have in particular a contraction cell in $PA$
	\[[\alpha \circ \beta, \beta \circ \alpha] : \alpha \circ \beta \lra \beta \circ \alpha\]
It is routine to check the axioms for a braided monoidal category using the contraction conditions at the top dimension; we show further that the symmetry axiom must hold, that is:
	\[\gamma_{\beta,\alpha} \circ \gamma_{\alpha, \beta} = 1\]
i.e.
	\[[\beta \circ \alpha, \alpha \circ \beta] \circ [\alpha \circ \beta, \beta \circ \alpha] = 1.\]
This is also true by contraction; in fact for any contraction 3-cell $[x,y]$ we have
	\[[x,y] \circ [y,x] = 1\]
since
	\[\phi([x,y] \circ [y,x]) = \phi[x,y] \circ \phi[y,x] = 1 = \phi(1).\]
Thus we see that a doubly degenerate 3-category is forced to be a symmetric monoidal category, not just a braided monoidal category as originally expected.

\section{Comparison with braided monoidal categories}\label{three}

In this section we give a precise sense in which Penon 3-categories are not the same as classical tricategories.  We exhibit a tricategory which does not arise as a Penon 3-category.  We will later show that this problem does not arise in the non-reflexive version.

The tricategory we examine is a doubly degenerate one: the free braided (strict) monoidal category on one object.  We show that its underlying monoidal category cannot be equipped with a symmetry and thus that it cannot be expressed as a doubly degenerate Penon 3-category.

Let \cl{B} denote the free braided (strict) monoidal category on one object.  This has

\begin{itemize}
\item objects the natural numbers
\item homsets $\cl{B}(n,m) = \left\{\begin{array}{ll}
	\mbox{$n$th braid group} & \mbox{if $m=n$} \\
	\emptyset & \mbox{otherwise}
	\end{array}\right.$
\item tensor product on objects addition, on morphisms juxtaposition of braids
\item unit object 0
\item braiding $\gamma_{m,n}: m+n \lra n+m$ depicted by

\bigskip

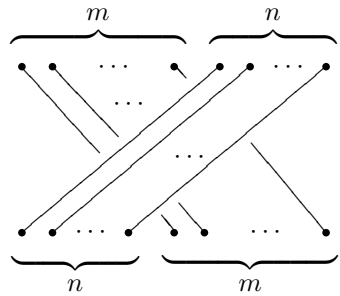
\[\xy
(10,8)*{\nowdot}="b1";
(14,8)*{\nowdot}="b2";
(19,8)*{\cdots};
(24,8)*{\nowdot}="b3";
(30,8)*{\nowdot}="b4";
(34,8)*{\nowdot}="b5";
(42,8)*{\cdots};
(50,8)*{\nowdot}="b6";
(10,30)*{\nowdot}="t1";
(14,30)*{\nowdot}="t2";
(22,30)*{\cdots};
(30,30)*{\nowdot}="t3";
(36,30)*{\nowdot}="t4";
(40,30)*{\nowdot}="t5";
(45,30)*{\cdots};
(50,30)*{\nowdot}="t6";
(20,32)*[@!270]{\left\{\begin{array}{l}
\vspace{20mm}
\end{array}\right.};
(20,37)*{m};
(43,32)*[@!270]{\left\{\begin{array}{l}
\vspace{14mm}
\end{array}\right.};
(43,37)*{n};
(17,6)*[@!90]{\left\{\begin{array}{l}
\vspace{14mm}
\end{array}\right.};
(17,1)*{n};
(40,6)*[@!90]{\left\{\begin{array}{l}
\vspace{20mm}
\end{array}\right.};
(40,1)*{m};
"t4"*{};"b1"*{} **\crv{}
	\POS?(.57)*{\hole}="x3"
	\POS?(.83)*{\hole}="x2"
	\POS?(.28)*{\hole}="x1";
"t5"*{};"b2"*{} **\crv{};
"t6"*{};"b3"*{} **\crv{}
	\POS?(.86)*{\hole}="y3"
	\POS?(.9)*{\hole}="y2"
	\POS?(.53)*{\hole}="y1";
"t1"*{};"x3" **\crv{};
"y3";"b4"*{} **\crv{};
"t2"*{};"x2" **\crv{};
"y2";"b5"*{} **\crv{};
"t3"*{};"x1" **\crv{};
"y1";"b6"*{} **\crv{};
(32,18)*{\cdots};
(24,25)*{\cdots};

\endxy\]

\end{itemize}

%

\noindent Now we observe that $\gamma$ is not a symmetry: for example in the case $m=n=1$ the composite $\gamma_{n,m} \circ \gamma_{m,n}$ is depicted by
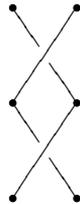
\[\xy
0;/r.25pc/:
(0,0)*{\nowdot}="b1";
(8,0)*{\nowdot}="b2";
(0,12)*{\nowdot}="m1";
(8,12)*{\nowdot}="m2";
(0,24)*{\nowdot}="t1";
(8,24)*{\nowdot}="t2";
"t2";"m1" **\crv{} \POS?(.5)*{\hole}="th";
"m2";"b1" **\crv{} \POS?(.5)*{\hole}="bh";
"t1";"th" **\crv{};
"th";"m2" **\crv{};
"m1";"bh" **\crv{};
"bh";"b2" **\crv{};
\endxy\]
which is not equal to the identity braid on 2.

\begin{prop}
The underlying monoidal category of \cl{B} cannot be equipped with a symmetry.
\end{prop}

\begin{prf}
We seek a symmetry
	\[\sigma_{AB}: A \otimes B \lra B \otimes A\]
natural in $A$ and $B$.  Put $A=B=1$.  In particular we need a morphism
	\[\sigma_{1,1}: 2 \lra 2.\]
The only such maps are given by $\gamma^{k}_{1,1}$ for all $k \in \bb{Z}$.  We have seen above that $\gamma_{1,1}$ is not a symmetry, and similarly $\gamma^{k}_{1,1}$ is not a braiding for any $k \geq 0$.  For $k=0$ we have the identity, but if $\sigma_{1,1} = \id$ then the braid axioms force $\sigma_{m,n} = \id$ for all $m, n$ which does not satisfy naturality.  Since there are no other morphisms $2 \lra 2 \in \cl{B}$ we conclude that there is no symmetry on this monoidal category. \end{prf}

%
%

We can now realise this braided monoidal category as a doubly degenerate tricategory whose 2-cells are the natural numbers and whose 3-cells $n \lra m$ are given by $\cl{B}(n,m)$ as above.  Composition along both 0-cells and 1-cells is given by $\otimes$ and the interchange constraint is derived from the braiding in the obvious way, with the homomorphism axiom following from the Yang-Baxter equation.

Note that this is in fact a Gray-category since everything is strict except interchange.  So we do indeed have a tricategory, and it does not satisfy the axioms for a Penon 3-category.

\subsubsection*{Remark}

We might ask if every tricategory is \emph{equivalent} to a Penon 3-category but this cannot be true.  The above braided monoidal category cannot be equivalent to a symmetric monoidal category since we know that a braided monoidal category is equivalent to a symmetric monoidal category if and only if it is itself symmetric \cite{gps1}.

\section{The non-reflexive case}\label{four}

In this section we give a non-reflexive version of Penon's definition.  This is a straightforward modification of the original definition, and is essentially the same as the version given in \cite{bat3}, but we present it slightly differently here.  We then show that the problems encountered in the previous sections no longer arise.

An $\omega$-category will now be defined as an algebra for a monad $N$ on ordinary (non-reflexive) globular sets.  We write a globular set $A$ as a diagram of sets as below.

\[\xy
(-20,10)*+{\cdots}="b4";
(0,10)*+{A(3)}="b3";
(20,10)*+{A(2)}="b2";
(40,10)*+{A(1)}="b1";
(60,10)*+{A(0)}="b0";
{\ar@<.5ex>^{s} "b4"; "b3"};
{\ar@<-.5ex>_{t} "b4"; "b3"};
{\ar@<.5ex>^{s} "b3"; "b2"};
{\ar@<-.5ex>_{t} "b3"; "b2"};
{\ar@<.5ex>^{s} "b2"; "b1"};
{\ar@<-.5ex>_{t} "b2"; "b1"};
{\ar@<.5ex>^{s} "b1"; "b0"};
{\ar@<-.5ex>_{t} "b1"; "b0"};

\endxy\]

Now we have a forgetful functor $\cl{Q} \lra \cat{GSet}$ given by composing the old forgetful functor $\cl{Q} \lra \cat{RGSet}$ with the forgetful functor $\cat{RGSet} \lra \cat{GSet}$.  As before, this has a left adjoint

\[\xy
(0,20)*++{\cl{Q}}="t0";
(0,0)*++{\cat{GSet}}="b0";
{\ar@<1.2ex>^{} "t0"; "b0"};
{\ar@<1.2ex>_{\adj}^{} "b0"; "t0"};
\endxy\]
inducing a monad which we will call $N$.  There is also an $n$-dimensional version as before, which we will call $N_n$.

\begin{Def} An $\omega$-category is an $N$-algebra.  An $n$-category is an $N_n$-algebra.  \end{Def}

Note that we thus have a commuting triangle of adjunctions

\[\xy
(0,0)*+{\cat{RGSet}}="bl";
(30,0)*+{\cat{GSet}}="br";
(15,18)*+{\cl{Q}}="t";
{\ar@<1.2ex>^{G_1}_{\top} "bl"; "br"};
{\ar@<1.2ex>^{F_1} "br"; "bl"};
(23,10)*[@!45]{\xy
(0,15)*++{}="t0";
(0,0)*++{}="b0";
{\ar@<1.2ex>^{} "t0"; "b0"};
{\ar@<1.2ex>_{\adj}^{} "b0"; "t0"};\endxy};
(20,7)*{{\scriptstyle F_3}};
(27,13)*{{\scriptstyle G_3}};
(7,10)*[@!315]{\xy
(0,15)*++{}="t0";
(0,0)*++{}="b0";
{\ar@<1.2ex>^{} "t0"; "b0"};
{\ar@<1.2ex>_{\adj}^{} "b0"; "t0"};\endxy};
(3,13)*{{\scriptstyle F_2}};
(10,7)*{{\scriptstyle G_2}};

\endxy\]

\noindent (where the bottom is monadic but the other two sides are not).  Thus we immediately have a construction of $NA$ -- we first add in putative identities freely and then proceed as in the reflexive case.

We may write $F_3A$ as

\[\left(
\def\objectstyle{\scriptstyle}
\def\labelstyle{\scriptstyle}
\vcenter{\xymatrix @-1.0pc  { NA \ar[d]^{} \\ T_RF_1A}} \right) =
\left(
\def\objectstyle{\scriptstyle}
\def\labelstyle{\scriptstyle}
\vcenter{\xymatrix @-1.0pc  { NA \ar[d]^{} \\ TA  }} \right)\]
where $T$ is the free strict $\omega$-category monad on non-reflexive globular sets, and we observe the crucial difference between the reflexive and non-reflexive versions of the theory --- in $NA$ the identities are freely added, so there are no non-identity cells which have identities as their source or target, whereas in $PA$ there may exist such cells.

\section{Degenerate 3-categories in the non-\\reflexive version}\label{five}

In this section we briefly examine degenerate 3-categories in the non-reflexive case and show that the previous problem of braidings being forced to be symmetries does not now arise.  We consider a doubly degenerate 3-category $\left(
\def\objectstyle{\scriptstyle}
\def\labelstyle{\scriptstyle}
\vcenter{\xymatrix @-1.0pc  { PA \ar[d]^{\theta} \\ A  }} \right)$.  As before, we construct a monoidal category from it with
	\begin{itemize}
	\item objects given by $A(2)$
	\item morphisms given by $A(3)$
	\item tensor product given by $\alpha \otimes \beta = \alpha \circ \beta$ as 2-cells of $A$.
	\end{itemize}
However, we cannot copy the previous construction of a braiding as we no longer have $\alpha \circ \beta = \beta \circ \alpha$ in the strict 3-category $TA$.  This is because $TA$ is now the free strict 3-category on the non-reflexive globular set $A$, so $TA(1) \neq 1$ although $A(1) = 1$.

In the reflexive version, the unique 1-cell of $A$ becomes the identity 1-cell of $T_RA$, so all the composites on the Eckmann-Hilton ``clockface" are equal.  In the non-reflexive version, the unique 1-cell of $A$ generates the 1-cells of $TA$ but there is a new (formal) 1-cell identity.  So the composites on the Eckmann-Hilton clockface do not even have the same source and target, and are certainly not equal.

This shows that the previous problem no longer arises; it remains to see how to construct a braiding at all.  We sketch a proposed argument here, but checking the axioms is not straightforward and we defer this to a future work.

Examining the Eckmann-Hilton clockface again we see that, apart from $\alpha \circ \beta$ and $\beta \circ \alpha$, the clock splits in two: the top half is all equal to $\beta * \alpha$ and the bottom half to  $\alpha * \beta$.  So in $NA$ we do have a contraction cell
	\[\chi = [(1 * \alpha) \circ (\beta * 1),\ (\beta * 1) \circ (1 * \alpha)]\]
(``10 o'clock to 2 o'clock"), so we seek to extend this to a braiding
	\[\alpha \circ \beta \lra  \beta \circ \alpha.\]
In the following argument we write $\circ$ and $*$ for the formal composition in $NA$, and evaluate these composites in $A$ by means of the algebra map $\theta$.  We write the unique 1-cell in $A$ as $e$ and the unit 1-cell in $NA$ as $I$.  Since this is only a sketch, we also ignore associativity issues with the understanding that for a precise construction these would need to be dealt with using further contraction cells.

We proceed in the following steps.

\numarabic
\begin{enumerate}

\item We have contraction 2-cells in $NA$ $\curlyl_e = [I*e,e]$, $\curlyr_e = [e*I,e]$ and also $[I*I, I]$.  We know that $\theta(I) = \theta(e) = e$, so by algebra associativity we have
	\[\theta(\curlyl_e) = \theta([I*I,I]) = \theta(\curlyr_e)\]
in $A$.  (This is the familiar result $\curlyl_I = \curlyr_I$ in any bicategory.)  By contraction, we also have pseudo-inverses for these cells, which we will denote $(\ )^*$.

\item We have contraction 3-cells in $NA$
	\[\lambda_\alpha = [\curlyl_e \circ (1_I \ast \alpha) \circ \curlyl_e^*, \alpha]\]
	\[\rho_\beta = [\curlyr_e \circ (\beta \ast 1_I) \circ \curlyr_e^*, \beta].\]

\item We now form $\rho \circ_1 \lambda$, composing these 3-cells along the 1-cell boundary, and apply $\theta$.  Now,
\[\xy
(12,7)*[@!90]{\left\{\begin{array}{l}
\vspace{18mm}
\end{array}\right.};
(0,10)*{s(\theta(\rho \circ_1 \lambda)) = \theta\left(\begin{array}{l}\vspace{3mm}\end{array}\right.\hspace{-5mm}
\theta(\curlyr_e) \circ \theta(\beta \ast 1_I) \circ \theta(\curlyr_e^*)
\circ \theta(\curlyl_e) \circ \theta(1_I \ast \alpha) \circ\theta(\curlyl_e^*)
\hspace{-5mm}\left.\begin{array}{l}\vspace{3mm}\end{array}\right)};
\endxy\]
so we can precompose by contraction cells at the middle factor (indicated), giving a composite 3-cell
\[\xymatrix{
\xi: \theta\left(\begin{array}{l}\vspace{3mm}\end{array}\right.\hspace{-5mm}\theta
(\curlyr_e)\circ\theta(\beta \ast 1_I) \circ \theta(1_I \ast \alpha) \circ \theta(\curlyl_e^*)\hspace{-5mm}\left.\begin{array}{l}\vspace{3mm}\end{array}
\right) \ar@3{->}[r] & \theta(\beta \circ \alpha)
}\]

\item Recall we have a contraction cell
\[\xymatrix{
\chi:\theta \bigleft (1_e \ast \alpha) \circ (\beta \ast 1_e) \bigright \ar@3{->}[r] & \theta\bigleft(\beta \ast 1_e) \circ (1_e \ast \alpha)\bigright
}.\]
Now, using algebra axioms we can rewrite this as
\[\xymatrix{
\theta \bigleft \theta(1_e \ast \alpha) \circ \theta(\beta \ast 1_e) \bigright \ar@3{->}[r] & \theta\bigleft \theta(\beta \ast 1_e) \circ \theta(1_e \ast \alpha)\bigright
}\]
and thus, to make it composable with $\xi$ it remains to compose it vertically with the identity 3-cells on $\curlyr_e$ and $\curlyl_e^*$; we have then bridged the ``gap" into 3 o'clock.

\item A similar argument then takes us from 9 o'clock to 10 o'clock.

\end{enumerate}

Evidently the above arguments are not ideal and we hope to find a more efficient method for calculating in this framework.  It remains to prove that this is in fact a braiding, but it is clear that the argument previously used to show that the braiding was a symmetry is no longer applicable.


\noindent Eugenia Cheng, Department of Pure Mathematics, University of Sheffield, Hicks Building, Hounsfield Road, Sheffield S3 7RH, UK.\\

\noindent Michael Makkai, Department of Mathematics, McGill University, Burnside Hall, Room 1005, 805 Sherbrooke Street West, Montreal, Quebec, Canada, H3A 2K6

\end{document}